%\nonstopmode
\magnification=1200

\def\og{\leavevmode\raise.3ex\hbox{$\scriptscriptstyle 
\langle\!\langle\,$}}
\def \fg {\leavevmode\raise.3ex\hbox{$\scriptscriptstyle 
\rangle\!\rangle\,\,$}}

\def\Z{{\bf Z}}

%%%%%%%%%%%%%%%%%%%%%%%%%%%%%%%%%%%%%%%%%%%%%%%%%%%%
% typoref.tex. V : January 18, 2000. 
% Author : Anthony PHAN
% Warning : syntaxe +- LaTeX 
% Sources :
% T. Lachand--Robert, ``La Ma\^\i trise de \TeX'',
% R\'ef\'erences crois\'ees;
% latex.ltx's sources;
% and of course the \TeX book.
%%%%%%%%%%%%%%%%%%%%%%%%%%%%%%%%%%%%%%%%%%%%%%%%%%%%%
%
\catcode`@=11
%
% style (look at the behavior of \item dans \bibitem too,
% and at one ,\  in \re@dreferenceslist)
% Feel free to change: 	\bibn@me (title like ``R\'ef\'erences'')
%			\bibliographym@rk (general style)
%
\def\bibn@me{R\'ef\'erences}
\def\bibliographym@rk{\centerline{{\sc\bibn@me}}
	\sectionmark\section{\ignorespaces}{\unskip\bibn@me}
	\bigbreak\bgroup
	\ifx\ninepoint\undefined\relax\else\ninepoint\fi}
%
% Beware of the \bgroup: it will be closed by \endthebibliography
%
% \refsp@ce is the spacing command that appens between multiple
% references.
%
\let\refsp@ce=\ 
\let\bibleftm@rk=[
\let\bibrightm@rk=]
%
% if you want more space between brackets...
%\let\refsp@ce=\thinspace
%\def\bibleftm@rk{[\thinspace}
%\def\bibrightm@rk{\thinspace]}
%
% frenchy stuff
%
\def\numero{n\raise.82ex\hbox{$\fam0\scriptscriptstyle o$}~\ignorespaces}
%
% new variables
%
\newcount\equationc@unt
\newcount\bibc@unt
\newif\ifref@changes\ref@changesfalse
\newif\ifpageref@changes\ref@changesfalse
\newif\ifbib@changes\bib@changesfalse
\newif\ifref@undefined\ref@undefinedfalse
\newif\ifpageref@undefined\ref@undefinedfalse
\newif\ifbib@undefined\bib@undefinedfalse
\newwrite\@auxout
%
% mark an equation
%
\def\eqnum{\global\advance\equationc@unt by 1%
\edef\lastref{\number\equationc@unt}%
\eqno{(\lastref)}}
%
% One can reference anything, just copy the former macro
% and use it so: \machin \label{truc}
% In machin you would have defined \lastref by some number
% or any text.
%
% References macros
%
% The next macros are the core of \ref and \cite commands.
% Its first argument may be ref, pageref or bib.
%
% It is too tricky to be explained.
% (It is a bit recursive.)
% It allows using \cite or \ref or ...
% with arbitrary many arguments,
% for instance:
% \cite{knuth1,knuth2,ma pomme}
%
% First argument is always ref, pageref or bib.
%
\def\re@dreferences#1#2{{%
	\re@dreferenceslist{#1}#2,\undefined\@@}}
\def\re@dreferenceslist#1#2,#3\@@{\def\next{#2}%
	\expandafter\ifx\csname#1@@\meaning\next\endcsname\relax
	??\immediate\write16
	{Warning, #1-reference "\next" on page \the\pageno\space
	is undefined.}%
	\global\csname#1@undefinedtrue\endcsname
	\else\csname#1@@\meaning\next\endcsname\fi
	\ifx#3\undefined\relax
	\else,\refsp@ce\re@dreferenceslist{#1}#3\@@\fi}
%
% notice that the former ``,\refsp@ce'' will separate
% multiple arguments. But beware of spaces
% while defining a reference or calling for it!
%
% tricky thing: \newlabel has two arguments
% {labelname}{{\lastref}{\pageref}}
% The second argument is read as two arguments
% by \newl@bel. This was necessary to get
% a jobname.aux containing the same syntax
% LaTeX would produce and use.
%
\def\newlabel#1#2{{\def\next{#1}\newl@bel#2}}
\def\newl@bel#1#2{%
	\expandafter\xdef\csname ref@@\meaning\next\endcsname{#1}%
	\expandafter\xdef\csname pageref@@\meaning\next\endcsname{#2}}
\def\label#1{{%
	\toks0={#1}\message{ref(\lastref) \the\toks0,}%
	\ignorespaces\immediate\write\@auxout%
	{\noexpand\newlabel{\the\toks0}{{\lastref}{\the\pageno}}}%
	\def\next{#1}%
	\expandafter\ifx\csname ref@@\meaning\next\endcsname\lastref%
	\else\global\ref@changestrue\fi%
	\newlabel{#1}{{\lastref}{\the\pageno}}}}
\def\ref#1{\re@dreferences{ref}{#1}}
\def\pageref#1{\re@dreferences{pageref}{#1}}
%
% bibliography macros
%
\def\bibcite#1#2{{\def\next{#1}%
	\expandafter\xdef\csname bib@@\meaning\next\endcsname{#2}}}
\def\cite#1{\bibleftm@rk\re@dreferences{bib}{#1}\bibrightm@rk}
%
% The argument of \beginthebibliography
% is any sequence of numerals which will represent
% the maximum \item's length. If you have less than 9
% \bibitem's, this argument may be {any numeral}.
% if you have between 100 and 999 \bibitem's
% this argument may be {any three numerals},
% and so on.
%
\def\beginthebibliography#1{\bibliographym@rk
	\setbox0\hbox{\bibleftm@rk#1\bibrightm@rk\enspace}
	\parindent=\wd0
	\global\bibc@unt=0
	\def\bibitem##1{\global\advance\bibc@unt by 1
		\edef\lastref{\number\bibc@unt}
		{\toks0={##1}
		\message{bib[\lastref] \the\toks0,}%
		\immediate\write\@auxout
		{\noexpand\bibcite{\the\toks0}{\lastref}}}
		\def\next{##1}%
		\expandafter\ifx
		\csname bib@@\meaning\next\endcsname\lastref
		\else\global\bib@changestrue\fi%
		\bibcite{##1}{\lastref}
		\medbreak
		\item{\hfill\bibleftm@rk\lastref\bibrightm@rk}%
		}
	}
\def\endthebibliography{\egroup\par}
%
% THE NEXT MACRO MUST BE INCLUDED
% IN THE \BYE COMMAND. FOR INSTANCE:
%
% \catcode`@=11
% \outer\def\bye{\@closeaux
% 	\par\vfill\supereject\end}
% \catcode`@=12
%
\def\@closeaux{\closeout\@auxout
	\ifref@changes\immediate\write16
	{Warning, changes in references.}\fi
	\ifpageref@changes\immediate\write16
	{Warning, changes in page references.}\fi
	\ifbib@changes\immediate\write16
	{Warning, changes in bibliography.}\fi
	\ifref@undefined\immediate\write16
	{Warning, references undefined.}\fi
	\ifpageref@undefined\immediate\write16
	{Warning, page references undefined.}\fi
	\ifbib@undefined\immediate\write16
	{Warning, citations undefined.}\fi}
%
% initialization of jobname.aux
%
\immediate\openin\@auxout=\jobname.aux
\ifeof\@auxout \immediate\write16
  {Creating file \jobname.aux}
\immediate\closein\@auxout
\immediate\openout\@auxout=\jobname.aux
\immediate\write\@auxout {\relax}%
\immediate\closeout\@auxout
\else\immediate\closein\@auxout\fi
%
% Let's read this file and open it out
%
\input\jobname.aux
\immediate\openout\@auxout=\jobname.aux
% this file will be closed by \bye.
%
% That's all, folks!
%
\catcode`@=12

\catcode`@=11
\def\bibliographym@rk{\bgroup}
%
% \bye est modifie pour la biblio et la table des matieres
%
\outer\def\bye{ 	\par\vfill\supereject\end}

\def\house#1{\setbox1=\hbox{$\,#1\,$}%
\dimen1=\ht1 \advance\dimen1 by 2pt \dimen2=\dp1 \advance\dimen2 by 2pt
\setbox1=\hbox{\vrule height\dimen1 depth\dimen2\box1\vrule}%
\setbox1=\vbox{\hrule\box1}%
\advance\dimen1 by .4pt \ht1=\dimen1
\advance\dimen2 by .4pt \dp1=\dimen2 \box1\relax}

  \def\eps{{\varepsilon}}

  \def\noi{\noindent}

\def\build#1_#2^#3{\mathrel{\mathop{\kern 0pt#1}\limits_{#2}^{#3}}}

\def\date {le\ {\the\day}\ \ifcase\month\or janvier
\or fevrier\or mars\or avril\or mai\or juin\or juillet\or
ao\^ut\or septembre\or octobre\or novembre
\or d\'ecembre\fi\ {\oldstyle\the\year}}

\font\fivegoth=eufm5 \font\sevengoth=eufm7 \font\tengoth=eufm10

\newfam\gothfam \scriptscriptfont\gothfam=\fivegoth
\textfont\gothfam=\tengoth \scriptfont\gothfam=\sevengoth

\def\smallsquare{\vbox{\hrule\hbox{\vrule height 1 ex\kern 1 ex\vrule}\hrule}}
\def\cqfd{\hfill \smallsquare\vskip 3mm}

%%%%%%%%%%%%%%%%%%%%%%%%%%%%%%%%%%%%%%%%%%%%

\centerline{}

\vskip 4mm

\centerline{
\bf Continued fractions and transcendental numbers}

\vskip 8mm
\centerline{Boris A{\sevenrm DAMCZEWSKI}, 
Yann B{\sevenrm UGEAUD}, and Les D{\sevenrm AVISON}}

\vskip 6mm

\vskip 8mm

\centerline{\bf 1. Introduction}

\vskip 6mm

It is widely believed that the 
continued fraction expansion of every irrational algebraic number 
$\alpha$ either is eventually periodic 
(and we know that this is the case if and only if $\alpha$ is a 
quadratic irrational), or it contains arbitrarily large partial quotients. 
Apparently, this question was first considered by 
Khintchine in  \cite{Khintchine} 
(see also \cite{Allouche,Shallit_survey,Waldschmidt} 
for surveys including a discussion on this subject). 
A preliminary step towards its resolution consists in providing explicit 
examples of transcendental continued fractions. 

The first result of this type goes back to  the 
pioneering work of 
Liouville \cite{Liouville}, who constructed transcendental 
real numbers with a very fast growing sequence of partial quotients. 
Indeed, the so-called `Liouville inequality' implies the 
transcendence of real numbers with very large partial quotients. 
Replacing it by Roth's theorem yields refined results, 
as shown by Davenport and Roth \cite{DavRoth}. In \cite{AdBuMB}, 
the argument of Davenport and Roth is slightly improved and Roth's 
theorem is replaced by a more recent result of Evertse \cite{Evertse97} 
to obtain the best known result of this type. Note that the constant $e$, 
whose continued fraction expansion is given by 
$$
e=[2;1,2,1,1,4,1,1,6,1,1,8,1,1,10,\ldots, 1, 1, 2n, 1, 1, \ldots],
$$ 
(see for instance \cite{Lang}) provides an explicit example of a transcendental number
with unbounded partial quotients; however, its transcendence does not
follow from the criteria from \cite{Liouville,DavRoth,AdBuMB}.

At the opposite side, there is a quest for finding 
explicit examples of transcendental continued fractions 
with bounded partial quotients. The first examples of such 
continued fractions were found by Maillet \cite{Maillet} 
(see also Section 34 of Perron \cite{Perron}). 
The proof of Maillet's results is based on a general form 
of the Liouville inequality 
which limits the approximation of algebraic numbers by quadratic irrationals. 
They were subsequently improved upon
by A. Baker \cite{Baker62,Baker64}, who used
the Roth theorem for number fields obtained by LeVeque \cite{Leveque56}.
Later on, Davison \cite{Dav89} applied a result of W. M. Schmidt \cite{Schm67}, 
saying that a real algebraic number cannot
be too well approximable by quadratic numbers, to show the transcendence of
some specific continued fractions (see Section 4). With the
same auxiliary tool, M. Queff\'elec \cite{Queffelec98} 
established the nice result that the
Thue--Morse continued fraction is transcendental (see Section 5). This
method has then been made more explicit, and combinatorial
transcendence criteria based on Davison's approach were given in
\cite{ADQZ,Dav,Bax03,LiSt}.

\medskip

Recently, Adamczewski and Bugeaud \cite{AdBuAM} obtained two new
combinatorial transcendence criteria for continued fractions 
that are recalled in Section 2. 
The main novelty in their approach is the use of a stronger Diophantine result
of W. M. Schmidt \cite{Schmidt72a,Schmidt80}, commonly known as 
the Subspace Theorem. After some work, it
yields considerable improvements upon the
criteria from \cite{ADQZ,Dav,Bax03,LiSt}. 
These allowed them to prove that the continued fraction 
expansion of every real algebraic number of degree at least three
cannot be `too simple', in various senses.
It is the purpose of the present work to give further
applications of their criteria to well-known (families of) continued
fractions. In addition, we provide a slight refinement of their
second criterion. 
Note  also that  a significative improvement of 
the results of Maillet and of Baker mentioned above is  
obtained in \cite{AdBuMB} thanks to a similar use of the Subspace Theorem.

\medskip

The present paper is organized as follows. In Section 2, 
we state the two main transcendence criteria of \cite{AdBuAM}, 
namely Theorems A and B. Then, Section 3 is devoted to a slight 
sharpening of Theorem B. In the subsequent Sections, we give various 
examples of applications of these transcendence criteria. 
We first solve in Section 4 a problem originally tackled by 
Davison in \cite{Dav89} and later considered by several 
authors in \cite{ADQZ,Dav,Bax03}. The
Rudin--Shapiro and the Baum--Sweet continued fractions are 
proved to be transcendental in Section 5. 
Then, Sections 6 and 7 are respectively devoted to folded continued
fractions and continued fractions arising from perturbed
symmetries (these last sequences were introduced 
by Mend\`es France \cite{Mendes}).  
In the last Section, we show through the study of another family of
continued fractions that, in some cases, rather than
applying roughly Theorem B, it is much better to go into its
proof and to evaluate continuants carefully.

\vskip 8mm

\centerline{\bf 2. The transcendence criteria}

\vskip 6mm

In this Section, we recall the transcendence criteria, 
namely Theorems A and B, obtained by the first two authors in \cite{AdBuAM}.

\medskip

Before stating the new criteria, we need to introduce some notation.
Let ${\cal A}$ be a given set, not necessarily finite. 
The length of a word
$W$ on the alphabet ${\cal A}$, that is, the number of letters
composing $W$, is denoted by $\vert W\vert$.
For any positive integer $k$, we write
$W^k$ for the word $W\ldots W$ ($k$ times repeated concatenation
of the word $W$). More generally, for any positive rational number
$x$, we denote by $W^x$ the word
$W^{[x]}W'$, where $W'$ is the prefix of
$W$ of length $\left\lceil(x- [x])\vert W\vert\right\rceil$. 
Here, and in all what follows, $[y]$ and
$\lceil y\rceil$ denote, respectively, the integer part and the upper
integer part of the real number $y$. 
Let ${\bf a}=(a_\ell)_{\ell \ge 1}$ 
be a sequence of elements from ${\cal A}$,
that we identify with the infinite word $a_1 a_2 \ldots a_\ell \ldots$.
Let $w$ be a rational number with $w>1$.
We say that ${\bf a}$ 
satisfies Condition $(*)_w$ if ${\bf a}$ is not
eventually periodic and if there exists 
a sequence of finite words $(V_n)_{n \ge 1}$ such that:

\medskip

\item{\rm (i)} For any $n \ge 1$, the word $V_n^w$ is a prefix
of the word ${\bf a}$;

\smallskip

\item{\rm (ii)} The sequence $(\vert V_n\vert)_{n \ge 1}$ is 
increasing.

\medskip

Roughly speaking, ${\bf a}$ satisfies Condition $(*)_w$
if ${\bf a}$ is not eventually periodic and if there exist infinitely 
many `non-trivial' repetitions (the size of which is measured
by $w$) at the beginning of the infinite word $a_1 a_2 \ldots a_\ell \ldots$

A first transcendence criterion for `purely' stammering continued
fractions is given in \cite{AdBuAM}.

\proclaim Theorem A. 
Let ${\bf a}=(a_\ell)_{\ell \ge 1}$ be a sequence of positive integers.
Let $(p_\ell/q_\ell)_{\ell \ge 1}$ denote the sequence of convergents to 
the real number
$$
\alpha:= [0; a_1, a_2, \ldots, a_\ell,\ldots].
$$ 
If there exists a rational number $w \ge 2$ such that
${\bf a}$ satisfies Condition $(*)_w$, then $\alpha$ is transcendental.
If there exists a rational number $w>1$ such that
${\bf a}$ satisfies Condition $(*)_w$, and if
the sequence $(q_\ell^{1/\ell})_{\ell \ge 1}$ is bounded
(which is in particular the case when the sequence ${\bf a}$ is bounded),
then $\alpha$ is transcendental.

Unfortunately, in the statement of Theorem A,
the repetitions must appear at the very beginning of ${\bf a}$.
When this is not the case, but the repetitions occur not
too far from the beginning of ${\bf a}$, then we have another
criterion.

Keep the above notation.
Let $w$ and $w'$ be non-negative rational numbers with $w>1$.
We say that ${\bf a}$ 
satisfies Condition $(**)_{w, w'}$ if ${\bf a}$ is not
eventually periodic and if there exist 
two sequences of finite words $(U_n)_{n \ge 1}$,   
$(V_n)_{n \ge 1}$ such that:

\medskip

\item{\rm (i)} For any $n \ge 1$, the word $U_nV_n^w$ is a prefix
of the word ${\bf a}$;

\smallskip

\item{\rm (ii)} The sequence
$({\vert U_n\vert} / {\vert V_n\vert})_{n \ge 1}$ is bounded 
from above by $w'$;

\smallskip

\item{\rm (iii)} The sequence $(\vert V_n\vert)_{n \ge 1}$ is 
increasing.

\medskip

We are now ready to present a transcendence criterion for (general)
stammering continued fractions, as stated in \cite{AdBuAM}.

\proclaim Theorem B. 
Let ${\bf a}=(a_\ell)_{\ell \ge 1}$ be a sequence of positive integers.
Let $(p_\ell/q_\ell)_{\ell \ge 1}$ denote the sequence of convergents to 
the real number
$$
\alpha:= [0; a_1, a_2, \ldots,a_\ell,\ldots].
$$ 
Assume that the sequence $(q_\ell^{1/\ell})_{\ell \ge 1}$ is bounded
and set $M = \limsup_{\ell \to + \infty} \, q_{\ell}^{1/\ell}$ and
$m = \liminf_{\ell \to + \infty} \, q_{\ell}^{1/\ell}$.
Let $w$ and $w'$ be non-negative real numbers with
$$
w > w' \biggl( 2 {\log M \over \log m} - 1 \biggr)
+ {\log M \over \log m}. \eqno (2.1)
$$
If ${\bf a}$ satisfies Condition $(**)_{w, w'}$, 
then $\alpha$ is transcendental.

It turns out that condition (2.1) can be slightly
weakened by means of a careful consideration of continuants.
This is the purpose of Section 3.

\vskip 6mm

\centerline{\bf 3. A slight sharpening of Theorem B}

\vskip 8mm

In this Section, 
we establish the following improvement of Theorem B.

\proclaim Theorem 3.1. 
Let ${\bf a}=(a_\ell)_{\ell \ge 1}$ be a sequence of positive integers.
Let $(p_\ell/q_\ell)_{\ell \ge 1}$ denote the sequence of convergents to 
the real number
$$
\alpha:= [0; a_1, a_2, \ldots,a_\ell,\ldots].
$$ 
Assume that the sequence $(q_\ell^{1/\ell})_{\ell \ge 1}$ is bounded
and set $M = \limsup_{\ell \to + \infty} \, q_{\ell}^{1/\ell}$ and
$m = \liminf_{\ell \to + \infty} \, q_{\ell}^{1/\ell}$.
Let $w$ and $w'$ be non-negative real numbers with
$$
w > w' \biggl( 2 {\log M \over \log m} - 1 \biggr) + 1. \eqno (3.1)
$$
If ${\bf a}$ satisfies Condition $(**)_{w, w'}$, 
then $\alpha$ is transcendental.

We note that the right hand side of (3.1) is always smaller than or 
equal to the right hand side of (2.1).

Before proceeding with the proof, we recall some useful facts
on continuants. For positive integers $a_1, \ldots, a_m$, denote 
by $K_m (a_1, \ldots, a_m)$ the denominator of the rational number
$[0; a_1, \ldots, a_m]$. It is commonly called a {\it continuant}.

\medskip

\proclaim Lemma 3.2.
For any positive integers $a_1, \ldots, a_m$ and any integer $k$ with
$1 \le k \le m-1$, we have
$$
K_m (a_1, \ldots , a_m) = K_m (a_m, \ldots, a_1)
$$
and
$$
\eqalign{
K_k (a_1, \ldots, a_k) \cdot K_{m-k} (a_{k+1}, \ldots, a_m)
& \le K_m (a_1, \ldots , a_m) \cr
& \le 2 \, K_k (a_1, \ldots, a_k) \cdot K_{m-k} (a_{k+1}, \ldots, a_m). \cr}
$$

\bigskip

\noindent {\bf Proof of Theorem 3.1.}
We follow step by step the proof of Theorem 2 from \cite{AdBuAM},
with a single modification. We keep the notation from \cite{AdBuAM}.
Recall that, for any $n \ge 1$, we set
$r_n = |U_n|$ and $s_n = |V_n|$.

Let $n$ be a positive integer. Let $\delta > 0$ be a (small)
real number. Since $w>1$ and $r_n \le w' s_n$, we get
$$
{2 r_n + \delta s_n \over r_n + (w-1) s_n} \le {2 w' s_n  + \delta s_n
\over w' s_n + (w-1) s_n} 
= {2 w' + \delta \over w' + w - 1} < {\log m \over \log M},
$$
by (3.1), if $\delta$ is sufficiently small. Consequently, there
exist positive real numbers $\eta$ and $\eta'$ with $\eta < 1$ such that
$$
2 (1 +\eta)(1 + \eta')r_n + \eta (1 + \eta') s_n <
(1 - \eta') (r_n + w s_n - s_n) {\log m \over \log M},  \eqno (3.2)
$$
for any $n \ge 1$.

In the course of the proof of Theorem 4.1
from \cite{AdBuAM} we had to bound from
above the quantity $q_{r_n}^{1 + \eta} \, q_{r_n + s_n}^{1 + \eta} 
\, q_{r_n + [w s_n]}^{-1}$.
Our new observation is that the estimate 
$$
q_{r_n + s_n} \, q_{r_n + [(w-1) s_n]} \ll
q_{r_n} \, q_{r_n + [w s_n]} 
$$
follows from Lemma 3.2 (here and below, the
numerical constant implied by $\ll$ does not depend on $n$). 
Consequently, we get
$$
q_{r_n}^{1 + \eta} \, q_{r_n + s_n}^{1 + \eta} \, q_{r_n + [w s_n]}^{-1}
\ll
q_{r_n}^{2 + \eta} \, q_{r_n + s_n}^{\eta} \, q_{r_n + [(w-1) s_n]}^{-1}.
\eqno (3.3)
$$

Assuming $n$ sufficiently large, we then have
$$
q_{r_n} \le M^{(1 + \eta') r_n}, \quad
q_{r_n + s_n} \le M^{(1 + \eta') (r_n+s_n)}, \quad {\rm and}
\quad q_{r_n + [(w-1) s_n]} \ge m^{(1 - \eta') (r_n + w s_n - s_n)},
$$
with $\eta'$ as in (3.2). Consequently, we get
$$
q_{r_n}^{2 + \eta} \, q_{r_n + s_n}^{\eta} \, q_{r_n + [(w-1) s_n]}^{-1}
\le M^{2(1 +\eta)(1 + \eta')r_n + \eta (1 + \eta') s_n} 
\, m^{- (1 - \eta') (r_n + w s_n - s_n)}
\le 1,
$$
by our choice of $\eta$ and $\eta'$. It then follows from (3.2) and (3.3) that
$$
q_{r_n}^{1 + \eta} \, q_{r_n + s_n}^{1 + \eta} \, q_{r_n + [w s_n]}^{-1} \ll 1,
$$
and, with the notation from \cite{AdBuAM}, we get the upper bound
$$
\prod_{1 \le j \le 4} \, |L_j ({\underline z_n})|
\ll ( q_{r_n} \, q_{r_n + s_n})^{-\eta}
$$
for any positive integer $n$. We then conclude as in that paper. \cqfd

\vskip 6mm

\centerline{\bf 4. Davison's continued fractions}

\vskip 8mm

Let $\theta$ be an irrational number with $0 < \theta < 1$.
Under some mild assumptions,
Davison \cite{Dav89} established the transcendence of the real number
$\alpha_{\theta}=[0;d_1,d_2,\ldots]$, with $d_n=1+([n\theta] \hbox{ mod }
2)$ for any $n \ge 1$. These extra assumptions were subsequently removed in 
\cite{ADQZ}. Then, Davison \cite{Dav} and Baxa \cite{Bax03}
studied the more general question of the transcendence of the real number
$\alpha_{k, \theta}=[0;d_1,d_2,\ldots]$, where $k \ge 2$
is an integer and $d_n=1+([n\theta] \hbox{ mod }
k)$ for any $n\ge 1$. The two authors obtained  
some partial results but their methods did not allow them to cover all the
cases. It turns out that Theorem A yields a complete
answer to this question. 

\proclaim Theorem 4.1. Let $\theta$ be an irrational number with
$0<\theta<1$ and let $k$ be an integer at
least equal to $2$. Let ${\bf d}=(d_n)_{n\geq 1}$ be defined by 
$d_n=1+([n\theta] \hbox{ mod }k)$ for any $n \ge 1$. Then, the number 
$\alpha_{k, \theta}=[0;d_1,d_2,\ldots]$ is transcendental.

In order to prove Theorem 4.1,
we need two auxiliary results (Lemmas 4.3 and 4.4 below),
which will be deduced from the following proposition obtained
in \cite{Dav}. Throughout this Section, $\theta
= [0; a_1, a_2, \ldots, a_n, \ldots]$ denotes an
irrational number in $(0, 1)$ and $(p_{n}/ q_{n})_{n \ge 0}$ is
the sequence of its convergents.

\proclaim Proposition 4.2. For any non-negative integers $n$ and $r$
with $1\leq r\leq q_{n+1}-1$, we have
$$
[(q_n+r)\theta]=p_n +[r\theta].
$$

Our first auxiliary result is the following.

\proclaim Lemma 4.3. For any integers $n$, $r$ and $s$
with $n\geq 1$, $0\leq s\leq a_{n+1}$ and $1\leq
r\leq q_n+q_{n-1}-1$, we have
$$
[(sq_n+r)\theta]=sp_n+[r\theta].
$$ 

\noi{\bf Proof.} The proof goes by induction on $s$. If $s=0$, the
statement is a tautology. Let us assume that the result holds for a
given $s$ with $0\leq s<a_{n+1}$. Then, we have $s+1\leq a_{n+1}$ and 
$$
(s+1)q_n+r\leq a_{n+1}q_n+q_n+q_{n-1}-1=q_{n+1}+q_n-1.
$$ 
Hence, $sq_n+r\leq q_{n+1}-1$ and Proposition 4.2 implies that
$$
[((s+1)q_n+r)\theta]=[(q_n+(sq_n+r))\theta]
=p_n+[(sq_n+r)\theta].
$$
By our inductive assumption, we thus obtain 
$$
[((s+1)q_n+r)\theta]=p_n+sp_n+[r\theta]=(s+1)p_n+[r\theta],
$$  
concluding the proof of the lemma.
\cqfd

\proclaim Lemma 4.4. For any integers $n$, $\ell$ and $r$
with $n\geq 1$, $\ell \geq 0$ and $1\leq
r\leq q_{n+1}-1$, we have
$$
[(q_{n+\ell}+q_{n+\ell-1}+\ldots+q_n+r)\theta]
= p_{n+\ell}+p_{n+\ell-1}+\ldots+p_n+[r\theta].
$$ 

\noi{\bf Proof.} The proof goes by induction on $\ell$. The case $\ell=0$ is
given by Proposition 4.2. Let us assume that the desired result holds for a
given non-negative integer $\ell$.
Let $n$ and $r$ be as in the statement of the lemma. 
Since $q_{m+2} \ge q_{m+1} + q_m$ holds
for any non-negative integer $m$, we easily get that
$$
q_{n+\ell}+\ldots+q_{n+1}+q_n+r\leq q_{n+\ell+2}-1.
$$
Then, we infer from Proposition 4.2 that
$$
\eqalign{
[(q_{n+\ell+1}+q_{n+\ell}+\ldots+q_n+r)\theta] & 
= [(q_{n+\ell+1}+(q_{n+\ell}+\ldots+q_n+r))\theta] \cr
& =p_{n+\ell+1}+[(q_{n+\ell}+\ldots+q_n+r))\theta]. \cr}
$$
Using the inductive assumption, this shows that
the desired result holds for $\ell+1$
and proves the lemma.  \cqfd 

{\noi{\bf Proof of Theorem 4.1.} We are now ready to prove Theorem
4.1 by showing that the sequence $(d_n)_{n\geq 1}$ satisfies 
Condition $(*)_w$ for a suitable real number $w > 1$. 
We first remark that the sequence ${\bf d}$ is not
eventually periodic since $\theta$ is irrational. 
We have to distinguish two cases. 

First, let us assume that there are
infinitely many integers $n$ such that $a_{n+1}\geq k$. 
For such an $n$, we infer from Lemma 4.3
with $s=k$ that
$$
[(kq_n+r)\theta] \equiv [r\theta] \hbox{ mod }k, 
\hbox{ for }1\leq r\leq q_n+q_{n-1}-1. \eqno (4.1)
$$  
Set $V_n=d_1\ldots d_{kq_n}$ and 
$W_n=d_{kq_n+1}d_{kq_n+2}\ldots
d_{(k+1)q_n+q_{n-1}-1}$ and view $V_n$ and $W_n$ as
words on the alphabet $\{1, 2, \ldots, k\}$. 
It follows from $(4.1)$ that $W_n$ is a prefix
of $V_n$. Furthermore, since 
$$
{\vert W_n\vert \over \vert V_n\vert}
={q_n+q_{n-1}-1\over kq_n}\geq{1\over k},
$$
the infinite word ${\bf d}$ begins in $V_n^{1+1/k}$.
Thus, ${\bf d}$ satisfies Condition $(*)_{1 + 1/k}$. 

Now, let us assume that there exists an integer $n_0$ such that
$a_n <  k$ for all $n\geq n_0$. Let $n \ge n_0$ be an integer.
At least two among the $k+1$ integers 
$p_n,p_n+p_{n+1},\ldots,p_n+p_{n+1}+\ldots+p_{n+k}$ are
congruent modulo $k$. Consequently, there exist 
integers $n'\geq n+1$ and $\ell$, with
$0 \le \ell \le k-1$ and $n' + \ell \le n + k$, such that 
$$
p_{n'+\ell}+p_{n'+\ell-1}+\ldots+p_{n'} \equiv 0\hbox{ mod }k.
$$
Set $N=q_{n'+\ell}+q_{n'+\ell-1}+\ldots+q_{n'}$. 
Then, Lemma 4.4 implies that
$$
[(N+r)\theta] \equiv [r\theta] \quad \hbox{ mod }k, \eqno (4.2)
$$
for $1\leq r\leq q_{n'+1}-1$. 
Set $V_n=d_1\ldots d_{N}$ and $W_n=d_{N+1}\ldots d_{N+q_{n'+1}-1}$
and view $V_n$ and $W_n$ as words on the alphabet $\{1, 2, \ldots, k\}$. 
Then, $(4.2)$ implies that $W_n$ is a prefix
of $V_n$. Since $q_{n'+\ell}\leq k^{\ell} q_{n'}$ (this follows
from the assumption $a_m < k$ for any $m \ge n_0$), we obtain
$$
{\vert W_n\vert \over \vert V_n\vert}
={q_{n'+1}-1\over N}\geq{q_{n'}\over kq_{n'+\ell}} \ge 
{1 \over k^{\ell+1}}\geq {1 \over k^{k}}.
$$
It thus follows that the infinite word ${\bf d}$ begins in
$V_n^{1+ 1/k^k}$. 

Consequently, the sequence ${\bf d}$ satisfies
in both cases Condition $(*)_{1+1/k^k}$. By Theorem A, the real number
$\alpha_{k, \theta}$ is transcendental, as claimed. \cqfd

\vskip 8mm

\centerline{\bf 5. Automatic continued fractions}

\vskip 6mm

Let $k \ge 2$ be an integer.
An infinite sequence ${\bf a}=(a_n)_{n\geq 0}$ is 
said to be $k$-automatic if $a_n$ is a finite-state function of the base-$k$ 
representation of $n$. This means that there exists
a finite automaton starting 
with the $k$-ary expansion of $n$ as input and producing the term $a_n$ as 
output. Finite automata are one of the most basic models of computation and   
take place at the bottom 
of the hierarchy of Turing machines. A nice reference on this topic is the 
book of Allouche and Shallit \cite{Allouche_Shallit}. 
We refer the reader to it for more details about this notion. 

Motivated by the Hartmanis-Stearns problem \cite{Hartmanis_Stearns}, 
the following question  
was addressed in \cite{AdBuAM}: do there exist algebraic numbers 
of degree at least three 
whose continued fraction expansion can be produced by a finite automaton?

The first result towards this problem is due to 
Queff\'elec \cite{Queffelec98}, who proved the transcendence 
of the Thue--Morse continued fractions. A more general statement 
can be found  in \cite{Queffelec00}. 
As it is shown in \cite{AdBuAM}, the transcendence of a large class 
of automatic continued fractions can be derived from Theorem A and B. 
In the present Section, we show how our transcendence
criteria apply to two others emblematic automatic sequences: 
the Rudin--Shapiro and the Baum--Sweet sequences.

\medskip

Before proving such results, we recall some classical facts about morphisms and morphic sequences. 
For a finite set ${\cal A}$, we denote by ${\cal A}^*$ the free monoid 
generated by ${\cal A}$. The empty word $\varepsilon$ is the neutral element 
of ${\cal A}^*$. Let ${\cal A}$ and ${\cal B}$ be two finite sets. An 
application from ${\cal A}$ to ${\cal B}^*$ can be uniquely extended to a 
homomorphism between the free monoids ${\cal A}^*$ and ${\cal B}^*$. We 
call morphism from ${\cal A}$ to ${\cal B}$ such a 
homomorphism. If there is a positive integer
$k$ such that each element of ${\cal A}$ is mapped to a word of
length $k$, then the morphism is called $k$-uniform or simply
uniform. 

A morphism 
$\sigma$ from ${\cal A}$ into itself is said to 
be prolongable if there exists a 
letter $a$ such that $\sigma(a)=aW$, where the word $W$ is such 
that $\sigma^n(W)$ is a non-empty word for every $n\geq 0$. 
In that case, the 
sequence of finite words $(\sigma^n(a))_{n\geq 1}$ converges in 
${\cal A}^{\Z_{\ge 0}}$ (endowed with the product topology of the discrete 
topology on each copy of ${\cal A}$) 
to an infinite word ${\bf a}$. Such an infinite word is clearly 
a fixed point for the map $\sigma$.

\vskip 6mm

{\bf 5.1 The Rudin--Shapiro continued fractions}

\vskip 6mm
  
Let ${\bf \varepsilon}=(\varepsilon_n)_{n\geq 0}$ be a
sequence with values in the set $\{+1,-1\}$. It is not difficult to
see that 
$$
\sup_{\theta\in [0,1]}\biggl \vert\sum_{0 \le n<N}
\varepsilon_ne^{2i\pi n\theta}\biggr \vert\geq\sqrt N,
$$
for any positive integer $N$.
In 1950, Salem asked the following question, related to some problems
in harmonic analysis: 
do there exist a sequence  ${\bf \varepsilon}=(\varepsilon_n)_{n\geq
0}$ in $\{+1,-1\}^{\Z_{\ge 0}}$ and a positive constant $c$ such that 
$$
\sup_{\theta\in [0,1]}\biggl \vert\sum_{0 \le n<N}
\varepsilon_ne^{2i\pi n\theta}\biggr \vert\leq c\sqrt N \eqno (5.1)
$$
holds for any positive integer $N$?
A positive answer to this problem was given by Shapiro
\cite{Shapiro} and Rudin \cite{Rudin}, who provided an explicit solution 
which is now known 
as the Rudin--Shapiro sequence. This sequence is 
a famous 
example of a $2$-automatic sequence and can be defined as follows: $r_n$ 
is equal to $+1$ (respectively $-1$) if the number of occurrences of
the pattern `$11$' in the binary representation of $n$ is even
(respectively odd). Theorem A yields the 
transcendence of the Rudin--Shapiro continued fractions.

\proclaim Theorem 5.1. Let $a$ and $b$ two distinct 
positive integers, and let ${\bf r}=(r_n)_{n\geq 0}$ be the
Rudin--Shapiro  
sequence on the alphabet $\{a,b\}$ (that is the symbol 
$1$ is replaced by $a$ and the symbol $-1$ is replaced
by $b$ in the usual Rudin--Shapiro sequence). Then, the real number 
$\alpha=[0;r_0,r_1,r_2,\ldots]$
is transcendental.

{\noi{\bf Proof.} We first infer from (5.1) that the Rudin--Shapiro 
sequence is not eventually periodic. 
We present now a useful description of this sequence.   
Let $\sigma$ be a morphism defined from $\{1,2,3,4\}^*$ into itself by: 
$\sigma(1)=12$, $\sigma(2)=42$, $\sigma(3)=13$ and
$\sigma(4)=43$. Let 
$$
{\bf u}=1242434213 \ldots
$$ 
be the fixed point of $\sigma$ begining with $1$ and let $\varphi$ 
be the morphism defined from $\{1,2,3,4\}^*$ to $\{a,b\}^*$ by: 
$\varphi(1)=\varphi(2)=a$ and $\varphi(3)=\varphi(4)=b$. 
It is known (see for instance \cite{Fog}, Ch. 5) 
that  
$$
{\bf r}=\varphi({\bf u}). \eqno (5.2)
$$ 
Since $\sigma^5 (1) = 1242434213$, ${\bf u}$ begins with $V^{1+1/8}$, 
where $V=12424342$. Then, it follows from (5.2) 
that for any positive integer $n$, the word ${\bf r}$ begins 
with $\varphi(\sigma^n(V^{1+1/8}))$. The morphism $\sigma$ 
being a uniform morphism, we easily check that 
$$
\varphi(\sigma^n(V^{1+1/8}))=(\varphi(\sigma^n(V)))^{1+1/8}.
$$ 
The sequence ${\bf r}$ thus satisfies the condition $(*)_{1+1/8}$. 
This ends the proof thanks to Theorem~A. \cqfd

\vskip 6mm

{\bf 5.2 The Baum--Sweet continued fractions}

\vskip 6mm

In 1976, Baum and Sweet \cite{Baum_Sweet1} proved that, unlike 
what is expected in the real case, the function field  
${\bf F}_2((X^{-1}))$ contains a cubic element (over ${\bf F}_2(X)$) 
with bounded partial quotient in its continued fraction expansion. 
This element is
$\sum_{n\ge 0} s_n X^{-n}$ , where $s_n$ 
is equal to $0$ if the binary representation of $n$ 
contains at least one string of $0'$s of odd length 
and $s_n$ is equal to $1$ 
otherwise. The sequence ${\bf s}=(s_n)_{n\geq 0}$ is now usually referred 
to as the Baum--Sweet sequence. 
It follows from \cite{ABL,AdBuAoM} that, 
for any integer $b\ge 2$, the real number $\sum_{n\ge 0}s_n/b^n$ 
is transcendental. 
Here, we use Theorem 3.1 to prove a similar result 
for the continued fraction expansion.

\proclaim Theorem 5.2. Let $a$ and $b$ be distinct 
positive integers, and let ${\bf s}=(s_n)_{n\geq 0}$ be the
Baum--Sweet sequence on the alphabet $\{a,b\}$ 
(that is, the symbol  $0$ is replaced by $a$ and the symbol $1$ 
is replaced by $b$ in the usual Baum--Sweet sequence). 
Then, the real number 
$\alpha=[0;s_0,s_1,s_2,\ldots]$
is transcendental.

\noi{\bf Proof.}  Let us first remark that the sequence ${\bf s}$ 
is not eventually periodic. Indeed, as shown in \cite{Baum_Sweet1},
the formal power series $\sum_{n\ge 0} s_n X^{-n}$ is 
a cubic element over ${\bf F}_2 (X)$, thus, it is
not a rational function. 

Let us now recall a useful description of 
the Baum--Sweet sequence on the alphabet $\{a, b\}$. Let 
$\sigma$ be the morphism defined from $\{1,2,3,4\}^*$ into itself by: 
$\sigma(1)=12$, 
$\sigma(2)=32$, 
$\sigma(3)=24$ and $\sigma(4)=44$. Let also $\varphi$ 
be the morphism defined from $\{1,2,3,4\}^*$ 
to $\{a,b\}^*$ by:  $\varphi(1)=\varphi(2)=b$ and $\varphi(3)=\varphi(4)=a$. 
Let $${\bf u} = 123224323244 \ldots$$ 
 denote the fixed point of $\sigma$ begining with $1$. It is known (see for instance \cite{Allouche_Shallit}, Ch. 6) 
that  
$${\bf s}=\varphi({\bf u}). \eqno (5.3)$$ 
Observe that ${\bf u}$ begins in the word $UV^{3/2}$, 
where $U=1$ and $V=232243$. Since the morphism $\sigma$ is
uniform, it follows from $(5.3)$ that, for any positive
integer $n$, the word ${\bf s}$ begins with $U_nV_n^{3/2}$, where
$U_n=\varphi(\sigma^n(U))$ and $V_n=\varphi(\sigma^n(V))$. In
particular, we have $\vert U_n\vert/\vert V_n\vert=1/6$ and 
${\bf s}$ thus satisfies Condition $(**)_{3/2, 1/6}$.
We further observe that the frequency of $a$ in the word $U_n$
tends to $1$ as $n$ tends to infinity. By Theorem 5 from \cite{ADQZ},
this implies that the sequence $(q_{\ell})^{1/\ell}_{\ell \ge 1}$
converges, where
$q_\ell$ denotes the denominator of
the $\ell$-th convergent of $\alpha$,
for any positive integer $\ell$. 
Thus, with the notation of Theorem 3.1, we have $M=m$.
Since $3/2 > 1 + 1/6$,
we derive from Theorem 3.1 that $\alpha$ is transcendental.
This finishes the proof of our theorem. \cqfd
\bigskip

\vskip 6mm

\centerline{\bf 6. Folded continued fractions} 
 
\vskip 8mm

Numerous papers, including the survey \cite{DeMevdP}, 
are devoted to paperfolding sequences. 
In this Section we consider folded continued fractions 
and we prove that they are always 
transcendental. 

A sheet of paper can be folded in half lengthways in
two ways: right half over left (the positive way) or left half over
right (the negative way). After having
been folded an infinite number of times, the sheet of paper can be
unfolded to display an infinite sequence of creases formed by hills
and valleys. For convenience, we will denote by $+1$ the hills and by
$-1$ the valleys. The simplest choice, that is to fold always in the
positive way, gives the well-known regular paperfolding sequence over 
the alphabet $\{+1,-1\}$. More
generally, if ${\bf e}=(e_n)_{n \ge 0}$ is in $\{+1,-1\}^{\Z_{\ge 0}}$, 
the associated paperfolding
sequence on the alphabet $\{+1,-1\}$ is obtained accordingly to the
sequence ${\bf e}$ of folding instructions, that is, the $n$-th fold
is positive if $e_n=+1$ and it is negative otherwise.

Among the numerous studies concerned with paperfolding sequences, much
attention has been brought on the way (quite intriguing) they are
related to   
the continued fraction expansion of some formal power series. 
Indeed, it can be shown that for any sequence 
${\bf e}=(e_n)_{n\geq 0}\in\{+1,-1\}^{\Z_{\ge 0}}$ of folded instructions,  
the continued fraction expansion of the formal power series 
$X\sum_{n\geq 0}e_nX^{-2^n}$ can be deduced from the associated paperfolding
sequence ${\bf f}$. 
As a consequence, the authors of \cite{vdP_Shallit} precisely 
described the continued fraction expansion 
of the real number $\xi = 2\sum_{n\geq 0}e_n2^{-2^n}$. 
In particular, they proved that such an expansion is also closely 
related to the sequence ${\bf f}$ and called it a
`folded continued fraction'. 
However, we point out that these `folded continued fractions' 
are not the same as those we consider in Theorem 6.1. 
Real numbers such as $\xi$
were shown to be transcendental in 
\cite{LvdP77} thanks to the so-called Mahler method 
(the observation that Ridout's Theorem \cite{Ridout} also implies
that such numbers are transcendental was 
for instance done in \cite{Adam}). 
Consequently, we get the transcendence of a family of continued 
fractions whose shape arises from paperfolding sequences; 
a fact mentioned in \cite{vdP_Shallit}. The following result 
has the same flavour though it is obtained in a totally different way.
It deals with another family of continued 
fractions arising from paperfolding sequences.

\proclaim Theorem 6.1. Let $a$ and $b$ be two positive distinct
integers, $(e_n)_{n\geq 0}
\in\{+1,-1\}^{\Z_{\ge 0}}$ be a sequence of folding instructions
and let ${\bf f}=(f_n)_{n\geq 0}$ be the associated paperfolding
sequence over the alphabet $\{a,b\}$ (that is, the symbol  $+1$ 
is replaced by $a$ and the symbol $-1$ is replaced by $b$ 
in the usual paperfolding sequence associated with $(e_n)_{n\geq 0}$). 
Then, the number 
$\alpha=[0;f_0,f_1,f_2,\ldots]$ is transcendental.

{\noi{\bf Proof.}  Let ${\bf f}=(f_n)_{n\geq 0}$ be the 
paperfolding sequence on the alphabet $\{a,b\}$ associated with the
sequence $(e_n)_{n\geq 0}$ in $\{+1,-1\}^{\Z_{\ge 0}}$ of folding instructions. 
We first notice that ${\bf f}$ is not eventually periodic
since no paperfolding sequence is eventually periodic (see  
\cite{DeMevdP}). 

We present now another useful description of the sequence ${\bf f}$.  
Let ${\cal F}_i$, $i\in\{+1,-1\}$, be the map defined from
the set $\{+1,-1\}^*$ 
into itself by 
$$
{\cal F}_i:\;\;w\longmapsto wi-({\overline w }),
$$
where $\overline{ w_1w_2\ldots w_m} := w_m w_{m-1}\ldots w_1$ 
denotes the mirror image of the word $w_1 w_2 \ldots w_m$, and
$-(w_1w_2\ldots w_m):=w'_1w'_2\ldots w'_m$, with $w'_i=+1$
(resp. $w'_i=-1$) if $w_i=-1$ (resp. $w_i=+1$). 
Let $\varphi$ be the morphism defined from $\{+1,-1\}^*$ 
to $\{a,b\}^*$ by $\varphi(+1)=a$ 
and $\varphi(-1)=b$. One can easily verify 
(see for instance \cite{DeMevdP}) that
$$
{\bf f}=\lim_{n\to + \infty}\varphi({\cal F}_{e_0}{\cal F}_{e_1}\ldots{\cal
F}_{e_n}(\varepsilon)), \eqno (6.1)
$$
where $\varepsilon$ denotes the empty word.

Let $n \ge 2$ be an integer and set $V_n={\cal
F}_{e_2}\ldots 
{\cal F}_{e_n}(\varepsilon)$. By (6.1), the finite word 
$$
F_n=\varphi({\cal F}_{e_0}{\cal F}_{e_1}\ldots{\cal
F}_{e_n}(\varepsilon))=\varphi({\cal F}_{e_0}{\cal F}_{e_1}(V_n))
$$
is a prefix of ${\bf f}$. Moreover, we have 
$$
{\cal F}_{e_0}{\cal F}_{e_1}(V_n)={\cal F}_{e_0}(V_n e_1-\overline{
V_n})=(V_n e_1-\overline{
V_n})e_0-(\overline{ V_n e_1-\overline{
V_n}}).
$$
This gives 
$$
{\cal F}_{e_0}{\cal F}_{e_1}(V_n)=(V_n e_1-\overline{
V_n}) e_0 V_n-(\overline{V_n e_1}).
$$
In particular, the word 
$$
\varphi(V_n e_1(-\overline{ V_n}) e_0 V_n)
$$
is a prefix of ${\bf f}$. 
Moreover, we have 
$$
\varphi(V_n e_1(-\overline{ V_n}) e_0 V_n)=
(\varphi(V_ne_1(-\overline{ V_n}) e_0 ))^{w_n},
$$
where $w_n=1+\vert \varphi(V_n)\vert / (2\vert \varphi(V_n)\vert +2)$. 
Thus, ${\bf f}$ begins with the word  
$(\varphi(V_ne_1(-\overline{ V_n}) e_0 ))^{5/4}$. Since $|\varphi(V_n)|$ 
tends to infinity with $n$,
this proves
that ${\bf f}$ satisfies Condition $(*)_{5/4}$. Consequently,
Theorem A implies that the real number $\alpha = [0; f_0, f_1, f_2, \ldots]$
is transcendental. This ends the proof. \cqfd
%%%%%%%%%%%%%%%%%%%%%%%%%%%%%%%%%%%%%%%%%%%%%

\vskip 6mm

\centerline{\bf 7. Generalized perturbed symmetry systems}

\vskip 8mm

A perturbed symmetry is a map defined from ${\cal A}^*$ into itself 
by $S_X(W)=WX\overline W$, where
${\cal A}$ is a finite set, $W$ and $X$ are finite words
on ${\cal A}$, and $\overline W$ is
defined as in the previous Section. 
Since $S_X(W)$ begins in $W$, we can
iterate the map $S_X$ to obtain the infinite sequence 
$$
S^{\infty}_X(W)=WX\overline{W}XW\overline{X}\overline{W}\ldots
$$
Such sequences were introduced by Mend\`es France \cite{Mendes}, who
proved that they are periodic if  and only if the word $X$ is a
palindrome, that is, $X={\overline X}$. 

More recently, Allouche and Shallit \cite{Allouche_Shallit98} 
generalized this notion as follows. For convenience, for any finite word $W$,
we set $W^E := W$ and $W^R := \overline{W}$.   
Let $k$ be a positive integer. Let
$X_1,X_2,\ldots,X_{k}$ be finite words of the same length $s$, and
let $e_1,e_2,\ldots,e_{k}$ be in $\{E,R\}$. Then, the associated
generalized perturbed symmetry is defined by
$$
S(W)=W\prod_{i=1}^{k} X_i W^{e_i}.
$$ 
Again, since $S(W)$ begins in $W$, we get
an infinite word $S^{\infty}(W)$ by iterating $S$. 
In \cite{Allouche_Shallit98}, the authors proved
that $S^{\infty}(W)$ is $(k+1)$-automatic and gave a necessary and
sufficient condition for this sequence to be eventually periodic. 

We introduce now a generalization of this process. 
First, let us assume in the previous definition that the
$X_i$ are finite and possibly empty words,
without any condition on their length. 
We can similarly define a map $S$ that we again call a generalized perturbed
symmetry. If all the $X_i$ have the same length, as previously, 
we may then call this map a uniform generalized perturbed symmetry. 
We point out that 
the  definition we consider here is more general than 
the one given in \cite{Allouche_Shallit98}. 
A generalized perturbed symmetry system is defined as a $4$-tuple 
$({\cal A},W,{\cal S},(S_n)_{n\geq 0})$, where ${\cal A}$ is a finite
set, $W\in{\cal A}^*$, $W\not=\varepsilon$, 
${\cal S}$ is a finite set of generalized
perturbed symmetries and $(S_n)_{n\geq 0}\in{\cal S}^{\Z_{\ge 0}}$. 
Then, it is easily verified that the sequence of finite
words 
$$
S_nS_{n-1}\ldots S_0(W)
$$ converges to an infinite
sequence, which we call the sequence produced by the generalized
perturbed symmetry system $({\cal A},W,{\cal S},(S_n)_{n\geq 0})$. 
When ${\cal A}$ is a subset of $\Z_{\ge 1}$, we obtain an infinite sequence
of positive integers that we can view as the continued fraction
expansion of a real number.

\proclaim Theorem 7.1. Generalized perturbed symmetry systems generate
either quadratic or transcendental continued fractions.

We point out that the set of numbers generated by 
generalized perturbed symmetry systems is not countable.

\medskip

{\noi{\bf Proof.} For any integer $n\geq 0$, we consider a 
generalized perturbed  symmetry $S_n$ asso\-ciated with the parameters 
$X_{1,n},X_{2,n},\ldots,X_{k_n,n}$ and $e_{1,n},e_{2,n},\ldots,e_{k_n,n}$
in $\{E, R\}$. 
Let ${\bf u}$ be the sequence generated by the generalized perturbed
symmetry system $({\cal A},W,{\cal S},(S_n)_{n\geq 0})$. For every positive 
integer $n$ we define the finite word  
$$
W_n=S_{n-1}\ldots S_0(W).
$$
Then, $W_n$ is a prefix of ${\bf u}$ and the sequence $(\vert
W_n\vert)_{n\geq 1}$ tends to infinity. 
To prove Theorem 7.1, we will show that there exists a rational $w$ 
greater than $1$, such that ${\bf u}$ satisfies 
Condition $(*)_w$. In order to do this, we have to distinguish two
different cases. 

First, assume that for
infinitely many integers $n$, we have $e_{1,n}=E$. In this case, 
$S_n(W_n)$ begins in $W_nX_{1,n} W_n^{e_{1,n}}=W_n X_{1,n} W_n$, thus the
sequence ${\bf u}$ begins in $W_n X_{1,n} W_n$. Since $\vert W_n\vert$ tends 
to infinity and $\vert
X_{1,n} \vert$ lies in a finite set, we have that ${\bf u}$ begins in
$(W_n X_{1,n})^{3/2}$ as soon as $n$ is large enough. Thus, ${\bf u}$
satisfies Condition $(*)_{3/2}$. 

Now, assume that there exists an integer $n_0$ such that 
for all $n\geq n_0$, we have $e_{1,n}=R$. Let $n$ be an integer greater
than $n_0$. Then, $S_n(W_n)$ begins in $W_nX_{1,n}W_n^R$, where
$$
W_n=S_{n-1}(W_{n-1}) = W_{n-1}
\prod_{i=1}^{k_{n-1}} X_{i,n-1}W_{n-1}^{e_{i, n-1}}.
$$
This implies that $S_n (W_n)$ begins with
$$
\eqalign{
& \biggl(W_{n-1}\prod_{i=1}^{k_{n-1}}
X_{i, n-1}W_{n-1}^{e_{i, n-1}}\biggr)
X_{1,n} \biggl(W_{n-1}\prod_{i=1}^{k_{n-1}}
X_{i,n-1} W_{n-1}^{e_{i,n-1}}\biggr)^R \cr
& =\biggl(W_{n-1}\prod_{i=1}^{k_{n-1}}X_{i,n-1} W_{n-1}^{e_{i,n-1}}\biggr)
X_{1,n} 
\biggl(\, \prod_{i=2}^{k_{n-1}} X_{i,n-1} W_{n-1}^{e_{i, n-1}}\biggr)^R
\biggl(W_{n-1} X_{1, n-1} W_{n-1}^{e_{1, n-1}}\biggr)^R. \cr}
$$
Since $n>n_0$, we have $e_{1,n-1}=R$. Consequently, $S_n(W_n)$ begins in 
$$
W_{n-1}\biggl(\, \prod_{i=1}^{k_{n-1}}
X_{i,n-1} W_{n-1}^{e_{i,n-1}}\biggr) X_{1,n}
\biggl(\, \prod_{i=1}^{k_{n-1}}
X_{i,n-1} W_{n-1}^{e_{i,n-1}} \biggr)^R
W_{n-1} =: W_{n-1} Y_{n-1} W_{n-1}.
$$
Since both the $k_j$ and the $X_{i,j}$ lie in a finite set, and since
$\vert W_n\vert$ tends to infinity with $n$, we get
that $S_n(W_n)$, and thus ${\bf u}$, 
begins in $(W_{n-1} Y_{n-1})^{1+ 1 / (3k)}$ for $n$
large enough, where we have set
$k=\max\{k_n : n \ge 1\}$. This implies 
that ${\bf u}$ satifies Condition $(*)_{1 + 1/(3k)}$, 
and we conclude the proof by applying Theorem A. \cqfd
%%%%%%%%%%%%%%%%%%%%%%%%%%%%%%%%%%%%%%%%%%%%%%

\vskip 6mm

\centerline{\bf 8. Inside the proof of Theorem 3.1: 
an example of the use of continuants} 
 
\vskip 8mm

The purpose of this Section is to point out that, in some cases, rather than
applying roughly Theorem 3.1, it is much better to go into its
proof and to evaluate continuants carefully. 
In order to illustrate this idea, we introduce a new 
family of continued fractions.

\medskip

Throughout this Section, we use the following notation.
Let $k \ge 3$ be an odd integer. Let $b_1, \ldots , b_k$
be positive integers with $b_1 < \ldots < b_k$. We
consider words defined on the alphabet $\Sigma = \{b_1, \ldots , b_k\}$.
The character $b_j$ denotes either the letter $b_j$, or the integer $b_j$,
according to the context. For any finite word $W$ on $\Sigma$
and for $j=1, \ldots, k$, let denote by $|W|_{b_j}$ the number of
occurrences of the letter $b_j$ in $W$. Set
$$
\Sigma^+_{\hbox{\sevenrm equal}} = \{ W \in \Sigma^+ : 
|W|_{b_j} = |W|_{b_1} \quad \hbox{for $2 \le j \le k$} \}.
$$

Let $\lambda > 1$ be real and, for any $n \ge 1$, let $W_n$ be in 
$\Sigma^+_{\hbox{\sevenrm equal}}$ with $|W_{n+1}| > \lambda |W_n|$. 
Consider the
infinite word ${\bf a} = (a_{\ell})_{\ell \ge 1}$ defined by
$$
{\bf a} = W_1 W_2^2 W_3^2 \ldots W_n^2 \ldots
$$
Observe that ${\bf a}$ is either ultimately periodic, or it
satisfies Condition $(*)_{w, w'}$ with $w=2$ and
$w' = 2/(\lambda - 1)$. 

If ${\bf a}$ is not eventually periodic, then
Theorem 3.1 implies the transcendence of the real number 
$\alpha = [0; a_1, a_2, \ldots ]$, provided 
$\lambda$ is sufficiently large in terms of $b_k$.
It turns out that the method
of proof of Theorem 3.1 is flexible enough to yield 
in some cases much better results: we can get a
condition on $\lambda$ that does not depend on the values of the $b_j$'s.

\proclaim Theorem 8.1. Keep the preceding notation.
Assume furthermore that $|W_n|$ is odd for all
sufficiently large $n$ and that $\lambda > 3.26$. Then, the real number
$$
\alpha = [0; a_1, a_2, \ldots ]
$$
is either quadratic, or transcendental.

In order to establish Theorem 8.1, we need three lemmas.
We keep the notation from Section 3 of \cite{ADQZ}.
In particular, for an integer matrix $A$, we denote by $\rho(A)$
its spectral radius and by $||A||$ its $L^2$-norm.
Recall that $\rho(A) = ||A||$ when $A$ is symmetrical.
Our first auxiliary result is extracted from \cite{ADQZ}.
In all what follows, we set $\gamma = 0.885$.

\proclaim Lemma 8.2. If $A = \pmatrix{ a & 1 \cr 1 & 0 \cr}$
and $B = \pmatrix{ b & 1 \cr 1 & 0 \cr}$ where $a$ and $b$ 
are distinct positive integers, then we have
$\rho(AB) > \bigl( \rho(A) \rho(B) \bigr)^{\gamma}$.

For $j=1, \ldots, k$, set $B_j = \pmatrix{ b_j & 1 \cr 1 & 0 \cr}$.
Set also
$$
X = {1 \over k} \, \sum_{j=1}^k \, \log \rho (B_j).
$$
For any finite word $V$ on $\Sigma$, denote by $K(V)$ the
corresponding continuant.

\proclaim Lemma 8.3. If $V$ is in $\Sigma^+_{\hbox{\sevenrm equal}}$, 
then we have
$$
{1 \over |V|} \, \log K(V) \le X.
$$

\medskip

\noindent {\bf Proof.}  Let $V = d_1 d_2 \ldots d_m$ be a finite word
defined over the alphabet $\Sigma^+_{\hbox{\sevenrm equal}}$. Set
$p_{m-1}/q_{m-1} = [0; d_1, \ldots, d_{m-1}]$ and
$p_m/q_m = [0; d_1, \ldots, d_m]$.
Then, we have $K(V) = q_m$  and
$$
K(V) \le \biggl \Vert \pmatrix{ q_m & q_{m-1} \cr p_m & p_{m-1} \cr}
\biggr \Vert.
$$
Setting $h_j = |V|_{b_j}$ for $j=1, \ldots, k$, it follows from
the theory of continued fractions that
$$
K(V) \le ||B_1||^{h_1} \ldots || B_k||^{h_k}.
$$
Since the $B_j$'s are symmetrical and 
$h_1 = \ldots = h_k = h$, we have 
$$K(V) \le \rho(B_1)^{h_1} \ldots \rho(B_k)^{h_k}. $$
Hence, the proof. \cqfd

Our last auxiliary result is the following.

\proclaim Lemma 8.4. If $V$ is in $\Sigma^+_{\hbox{\sevenrm equal}}$
with $|V|$ odd, then we have
$$
{1 \over |V|} \, \log K(V) > \gamma X - {\log 4 \over |V|}.
$$

\medskip

\noindent {\bf Proof.} We use repeatedly a particular case
of Theorem 3.4 from \cite{Dav}.
It asserts that if $W$ is the product of an odd number $m$ 
of matrices $B_1, \ldots , B_k$, each of which occurring 
exactly $\ell$ times, then we have
$$
{\rm tr} (W) \ge \rho(B_1 B_k)^{\ell} \,
{\rm tr} (W'),
$$
where $W'$ is the product arising from $W$ by replacing the
matrices $B_1$ and $B_k$ by the identity matrix. As usual,
${\rm tr}(M)$ denotes the trace of the matrix $M$.

With the notation of the proof of Lemma 8.3, we then get that 
$$
\eqalign{
K(V) & \ge {1 \over 2} \, {\rm tr} 
\pmatrix{ q_m & q_{m-1} \cr p_m & p_{m-1}}
\ge {1 \over 2} \, \rho(B_1 B_k)^h \ldots 
\rho (B_{(k-1)/2} B_{(k+3)/2})^h \, {\rm tr} (B_{(k+1)/2}^h) \cr
& \ge {1 \over 4} \, \rho(B_1 B_k)^h \ldots 
\rho (B_{(k-1)/2} B_{(k+3)/2})^h \rho (B_{(k+1)/2}^h) \cr
& > {1 \over 4} \, \bigl( \rho(B_1) \ldots \rho (B_k) \bigr)^{\gamma h}, \cr}
$$
by Lemma 8.2. The lemma follows. \cqfd

\medskip

We have now all the tools needed to establish Theorem 8.1.

\medskip

\noindent {\bf Proof of Theorem 8.1.} 
For any $n \ge 2$, set
$$
U_n = W_1 W_2^2 \ldots W_{n-1}^2 \quad
{\rm and} \quad V_n = W_n.
$$
Clearly, ${\bf a}$ begins in $U_n V_n^2$. 
Denote by $K(U_n)$ and by $K(V_n)$ the continuants associated to
the words $U_n$ and $V_n$, respectively. In view of Lemma 3.2 and of (3.3),
the theorem is proved as soon as we establish that there exists
a positive real number $\eps$ such that
$$
\log K(V_n) > (1 + \eps) \log K(U_n),  \eqno (8.1)
$$
for any sufficiently large integer $n$.

To prove (8.1), we first infer from Lemmas 3.2  
and 8.3 that
$$
\log K(U_n) < \log K(V_1) + 2 \, \sum_{j=2}^{n-1} \,
\log K(V_j) + 2 n < 2 X \, \sum_{j=1}^{n-1} \, |V_j| + 2 n.
$$
Consequently, we get
$$
{1 \over |V_n|} \, \log K(U_n) < {2X \over \lambda - 1} + 
{2 n \over |V_n|}. \eqno (8.2)
$$
On the other hand, Lemma 8.4 gives us that
$$
{1 \over |V_n|} \, \log K(V_n) > \gamma X - {\log 4 \over |V_n|}. \eqno (8.3)
$$
We then infer from (8.2) and (8.3) that (8.1) is satisfied for some
positive $\eps$ as soon as we have $\gamma > 2/(\lambda - 1)$,
that is, $\lambda > 1 + 2 \, \gamma^{-1} \asymp 3.25 \ldots$
This completes the proof of the theorem. \cqfd

%%%%%%%%%%%%%%%%%%%%%%%%%%%%%%%%%%%%%%%%%%%%%%%%%%%%%%%%%%%%%%%%%%%%%%%%%%%

\vskip 12mm

\centerline{\bf References}

\vskip 7mm

\beginthebibliography{999}

\bibitem{Adam}
B. Adamczewski,
{\it Transcendance \og \`a la Liouville \fg de certain nombres r\'eels},
C. R. Acad. Sci. Paris {\bf 338} (2004), 511--514.

\bibitem{AdBuAoM}
B. Adamczewski \& Y. Bugeaud,
{\it On the complexity of algebraic numbers I.
Expansions in integer bases}.
Annals of Math. To appear.

\bibitem{AdBuAM}
B. Adamczewski \& Y. Bugeaud,
{\it On the complexity of algebraic numbers II. Continued
fractions}, Acta Math. To appear.

\bibitem{AdBuMB}
B. Adamczewski \& Y. Bugeaud,
{\it On the Maillet--Baker continued fractions}. Preprint.

\bibitem{ABL}
B. Adamczewski, Y. Bugeaud \&  F. Luca,
{\it Sur la complexit\'e des nombres alg\'ebriques},
C. R. Acad. Sci. Paris {\bf 339} (2004), 11--14.

\bibitem{Allouche}
J.-P. Allouche,
{\it Nouveaux r\'esultats de transcendance de r\'eels \`a 
d\'eveloppements non al\'eatoire},
Gaz. Math. {\bf 84} (2000), 19--34.

\bibitem{ADQZ}
J.-P. Allouche, J. L. Davison, M. Queff\'elec \& L. Q. Zamboni,
{\it Transcendence of Sturmian or morphic continued fractions},
J. Number Theory {\bf 91} (2001), 39--66.

\bibitem{Allouche_Shallit98}
J.-P. Allouche \& J. O. Shallit, 
{\it Generalized Pertured Symmetry}, 
Europ. J. Combinatorics {\bf 19} (1998), 401--411.

\bibitem{Allouche_Shallit}
J.-P. Allouche \& J. O. Shallit, 
Automatic Sequences: Theory, Applications, Generalizations, 
Cambridge University Press, Cambridge, 2003.

\bibitem{Baker62} 
A. Baker,
{\it Continued fractions of transcendental numbers},
Mathematika {\bf 9} (1962), 1--8.

\bibitem{Baker64} 
A. Baker,
{\it On Mahler's classification of transcendental numbers},
Acta Math. {\bf 111} (1964), 97--120.

\bibitem{Baum_Sweet1}
L. E. Baum \& M. M. Sweet,
{\it Continued fractions of algebraic power series in characteristic
$2$},
Annals of Math. {\bf 103} (1976), 593--610.

\bibitem{Bax03}
C. Baxa,
{\it Extremal values of continuants and transcendence of certain
continued fractions},
Adv. in Appl. Math. {\bf 32} (2004), 754--790.

\bibitem{DavRoth}
H. Davenport \& K. F. Roth, 
{\it Rational approximations to algebraic numbers}, 
Mathematika {\bf 2} (1955), 160--167.

\bibitem{Dav89}
J. L. Davison,
{\it A class of transcendental numbers with bounded partial quotients}.
In R. A. Mollin, ed., Number Theory and Applications, pp. 365--371, Kluwer
Academic Publishers, 1989.

\bibitem{Dav}
J. L. Davison,
{\it Continued fractions with bounded partial quotients},
Proc. Edinburgh Math. Soc. {\bf 45} (2002), 653--671.

\bibitem{DeMevdP}
F. M. Dekking, M. Mend\`es France \& A. J. van der Poorten, 
{\it Folds!}, Math. Intelligencer {\bf 4} (1982), 130--138, 173--181,
190--195. Erratum, {\bf 5} (1983), 5. 

\bibitem{Evertse97}
J.-H. Evertse,
{\it The number of algebraic numbers of given degree approximating
a given algebraic number}. In: Analytic number theory (Kyoto, 1996),
53--83, London Math. Soc. Lecture Note Ser. 247,
Cambridge Univ. Press, Cambridge, 1997.
 
 \bibitem{Fog}
N. Pytheas Fogg, 
Substitutions in Dynamics, Arithmetics and Combinatorics,
Lecture Notes in Mathematics 1794, Springer-Verlag, 2002.
 
 \bibitem{Hartmanis_Stearns}
 J. Hartmanis \& R. E. Stearns,
{\it On the computational complexity of algorithms},
Trans. Amer. Math. Soc. {\bf 117} (1965), 285--306.

\bibitem{Khintchine}
A. Ya. Khintchine,
Continued fractions, Gosudarstv. Izdat. Tehn.-Theor. Lit. 
Moscow-Leningrad, 2nd edition, 1949 (in Russian).

\bibitem{Lang}
S. Lang,
Introduction to Diophantine Approximations, Sprin\-ger-Verlag, 1995.

\bibitem{Leveque56}
W. J. LeVeque, 
Topics in number theory, Vol. II, Addison-Wesley, 1956.

\bibitem{LiSt}
P. Liardet \& P. Stambul,
{\it S\'eries de Engel et fractions continu\'ees},
J. Th\'eor. Nombres Bordeaux {\bf 12} (2000), 37--68.

\bibitem{Liouville} J.~Liouville, {\it Sur des classes 
tr\`es \'etendues de quantit\'es 
dont la valeur n'est ni alg\'ebri\-que, ni m\^eme r\'eductible 
\`a des irrationelles} 
{\it alg\'e\-bri\-ques}, 
C. R. Acad. Sci. Paris {\bf 18} (1844), 883--885; 910-911.

\bibitem{LvdP77}
J. H. Loxton \& A. J. van der Poorten,
{\it Arithmetic properties of certain functions in several variables
III}, Bull. Austral. Math. Soc. {\bf 16} (1977), 15--47.

\bibitem{Maillet}
E. Maillet,
Introduction \`a la th\'eorie des nombres transcendants et des propri\'et\'es
arithm\'etiques des fonctions, Gauthier-Villars, Paris, 1906.

\bibitem{Mendes}
M. Mend\`es France,
{\it Principe de la sym\'etrie perturb\'ee}. 
In: {S\'eminaire de Th\'eorie des
Nombres, Paris 1979-80}, M.-J. Bertin (\'ed.), Birkh\"auser, Boston, 1981,
pp. 77--98.

\bibitem{Perron}
O. Perron,
Die Lehre von den Ketterbr\"uchen.
Teubner, Leipzig, 1929.

\bibitem{vdP_Shallit}
A. J. van der Poorten \& J. O. Shallit,
{\it Folded continued fractions}, 
J. Number Theory {\bf 40} (1992), 237--250.

\bibitem{Queffelec98}
M. Queff\'elec,
{\it Transcendance des fractions continues de Thue--Morse},
J. Number Theory {\bf 73} (1998), 201--211.

\bibitem{Queffelec00}
M. Queff\'elec,
{\it Irrational number with automaton-generated continued fraction expansion},
In J.-M. Gambaudo, P. Hubert, P. Tisseur, and S. Vaienti, editors, 
{\it Dynamical Systems: From Crystal to Chaos}, World Scientific,
2000, 190--198.

\bibitem{Ridout}
D. Ridout,
{\it Rational approximations to algebraic numbers},
Mathematika {\bf 4} (1957), 125--131.

\bibitem{Rudin}
W. Rudin,
{\it Some theorems on Fourier coefficients},
Proc. Amer. Math. Soc. {\bf 10} (1959), 855--859.

\bibitem{Schm67}
 W. M. Schmidt,
{\it On simultaneous approximations of two algebraic numbers by rationals},
Acta Math. {\bf 119} (1967), 27--50.

\bibitem{Schmidt72a}
 W. M. Schmidt,
{\it Norm form equations},
Annals of Math. {\bf 96} (1972), 526--551.
  
\bibitem{Schmidt80}
W. M. Schmidt, 
{\it Diophantine approximation},
Lecture Notes in Mathematics 785, Springer, Berlin, 1980.

\bibitem{Shallit_survey}
J. O. Shallit, 
{\it Real numbers with bounded partial quotients}, 
Enseign. Math. {\bf 38} (1992), 151--187.

\bibitem{Shapiro}
H. S. Shapiro,
{\it Extremal problems for polynomials and power series},
Master's thesis, MIT, 1952.

\bibitem{Waldschmidt}
M. Waldschmidt,
{\it Un demi-si\`ecle de transcendance}. In:
Development of mathematics 1950--2000, pp. 1121--1186, 
Birkh\"auser, Basel, 2000.  

\endthebibliography

\vskip 6mm

\noindent Boris Adamczewski   \hfill{Yann Bugeaud}

\noindent   CNRS, Institut Camille Jordan  
\hfill{Universit\'e Louis Pasteur}

\noindent   Universit\'e Claude Bernard Lyon 1 
\hfill{U. F. R. de math\'ematiques}

\noindent   B\^at. Braconnier, 21 avenue Claude Bernard
 \hfill{7, rue Ren\'e Descartes}

\noindent   69622 VILLEURBANNE Cedex (FRANCE)   
\hfill{67084 STRASBOURG Cedex (FRANCE)}

\vskip2mm
 
\noindent {\tt Boris.Adamczewski@math.univ-lyon1.fr}
\hfill{{\tt bugeaud@math.u-strasbg.fr}}

\vskip 6mm

\centerline{Les Davison}

\centerline{Department of Mathematics and Computer Science}

\centerline{Laurentian University}

\centerline{Sudbury, Ontario}

\centerline{CANADA P3E 2C6}

\vskip2mm

\centerline{{\tt ldavison@cs.laurentian.ca}}

\bye